\documentclass{article}

\usepackage{latexsym}

\usepackage{amsmath, amssymb}

\title{THE SCATTERING MATRIX WITH RESPECT TO AN HERMITIAN MATRIX OF A GRAPH}

\author{Takashi KOMATSU \\
Department of Bioengineering, School of Engineering,\\ 
The University of Tokyo \\
Bunkyo, Tokyo, 113-8656, JAPAN \\ 
e-mail: komatsu@coi.t.u-tokyo.ac.jp \\ 
Norio KONNO \\
Department of Applied Mathematics, Faculty of Engineering, \\ 
Yokohama National University \\
Hodogaya, Yokohama 240-8501, JAPAN \\
e-mail: konno-norio-bt@ynu.ac.jp \\ 
Iwao SATO \\ 
Oyama National College of Technology \\
Oyama, Tochigi 323-0806, JAPAN \\ 
e-mail: isato@oyama-ct.ac.jp }

\begin{document}
 \maketitle

\clearpage

\begin{abstract}
Recently, Gnutzmann and Smilansky \cite{GS} presented a formula for 
the bond scattering matrix of a graph with respect to a Hermitian matrix.  
We present another proof for this Gnutzmann and Smilansky's formula 
by a technique used in the zeta function of a graph. 
Furthermore, we generalize Gnutzmann and Smilansky's formula to a regular covering 
of a graph. 
Finally, we define an $L$-fuction of a graph, and present a determinant expression. 
As a corollary, we express the generalization of Gnutzmann and Smilansky's formula 
to a regular covering of a graph by using its $L$-functions. 
\end{abstract}

\vspace{5mm}


\vspace{5mm}

The contact author for correspondence:

Iwao Sato 

Oyama National College of Technology, 
Oyama, Tochigi 323-0806, JAPAN 

Tel: +81-285-20-2176

Fax: +81-285-20-2880

E-mail: isato@oyama-ct.ac.jp

\clearpage

\section{Introduction} 

Ihara zeta functions of graphs started from Ihara zeta functions of regular 
graphs by Ihara \cite{Ihara}. 
Originally, Ihara presented $p$-adic Selberg zeta functions of 
discrete groups, and showed that its reciprocal is a explicit polynomial. 
Serre \cite{Serre} pointed out that the Ihara zeta function is the zeta function of 
the quotient $T/ \Gamma $ (a finite regular graph) of the one-dimensional 
Bruhat-Tits building $T$ (an infinite regular tree) associated with 
$GL(2, k_p)$. 

A zeta function of a regular graph $G$ associated with a unitary 
representation of the fundamental group of $G$ was developed by 
Sunada \cite{Sunada1, Sunada2}. 
Hashimoto \cite{Hashimoto} treated multivariable zeta functions of bipartite graphs. 
Bass \cite{Bass} generalized Ihara's result on the Ihara zeta function of 
a regular graph to an irregular graph, and showed that its reciprocal is 
again a polynomial. 
Various proofs of Bass' Theorem were given by 
Stark and Terras \cite{ST}, Foata and Zeilberger \cite{FZ}, Kotani and Sunada \cite{KS}. 
Sato \cite{Sato} defined the second weighted zeta function of a graph by using not an infinite 
product but a determinant.

The spectral determinant of the Laplacian on a quantum graph 
is closely related to the Ihara zeta function of a graph(see \cite{CDT, Des, GS, Smi}). 
Smilansky \cite{Smi} considered spectral zeta functions and trace formulas 
for (discrete) Laplacians on ordinary graphs, and expressed 
some determinant on the bond scattering matrix of a graph $G$ 
by using the characteristic polynomial of its Laplacian. 
Recently, Gnutzmann and Smilansky \cite{GS} presented a formula for 
the bond scattering matrix of a graph with respect to a Hermitian matrix.  

In this paper, we another proof for the Gnutzmann and Smilansky's formula 
on the bond scattering matrix of a graph with respect to a Hermitian matrix.  
by a technique used in the zeta function of a graph, and treat some related topics.
In Section 2, we review the Ihara zeta function and the bond scattering matrix of a graph $G$. 
In Section 3, we present another proof for the Gnutzmann and Smilansky's formula 
by a technique used in the zeta function of a graph. 
In Section 4, we we express a new zeta function of $G$ on the bond scattering matrix 
of $G$ with respect to a Hermitian matrix by using the Euler product. 
In Section 5, we generalize the Gnutzmann and Smilansky's formula to a regular covering 
of $G$. 
In Section 6, we define an $L$-fuction of $G$, and present its determinant expression. 
As a corollary, we express the generalization of the Gnutzmann and Smilansky's formula 
to a regular covering of $G$ by using its $L$-functions.

\section{The zeta functions and the bond scattering matrix of a graph}

Graphs treated here are finite.
Let $G=(V(G),E(G))$ be a connected graph (possibly multiple edges and loops) 
with the set $V(G)$ of vertices and the set $E(G)$ of unoriented edges $uv$ 
joining two vertices $u$ and $v$. 
For $uv \in E(G)$, an arc $(u,v)$ is the oriented edge from $u$ to $v$. 
Set $D(G)= \{ (u,v),(v,u) \mid uv \in E(G) \} $. 
For $b=(u,v) \in D(G)$, set $u=o(b)$ and $v=t(b)$. 
Furthermore, let $b^{-1} =(v,u)$ be the {\em inverse} of $b=(u,v)$. 

A {\em path $P$ of length $n$} in $G$ is a sequence 
$P=(b_1, \cdots ,b_n )$ of $n$ arcs such that $b_i \in D(G)$,
$t( b_i )=o( b_{i+1} )(1 \leq i \leq n-1)$, 
where indices are treated $mod \  n$. 
Set $ \mid P \mid =n$, $o(P)=o( b_1 )$ and $t(P)=t( b_n )$. 
Also, $P$ is called an {\em $(o(P),t(P))$-path}. 
We say that a path $P=(b_1, \cdots ,b_n )$ has a {\em backtracking} 
or {\em back-scatter} if $ b_{i+1} =b^{-1}_i $ for some $i(1 \leq i \leq n-1)$. 
A $(v, w)$-path is called a {\em $v$-cycle} 
(or {\em $v$-closed path}) if $v=w$. 
The {\em inverse cycle} of a cycle 
$C=( b_1, \cdots ,b_n )$ is the cycle 
$C^{-1} =( \hat{b}_n , \cdots ,\hat{b}_1 )$.

We introduce an equivalence relation between cycles. 
Two cycles $C_1 =(e_1, \cdots ,e_m )$ and 
$C_2 =(f_1, \cdots ,f_m )$ are called {\em equivalent} if there exists 
$k$ such that $f_j =e_{j+k} $ for all $j$. 
The inverse cycle of $C$ is in general not equivalent to $C$. 
Let $[C]$ be the equivalence class which contains a cycle $C$. 
Let $B^r$ be the cycle obtained by going $r$ times around a cycle $B$. 
Such a cycle is called a {\em power} of $B$. 
A cycle $C$ is {\em reduced} if both $C$ and $C^2 $ have no backtracking. 
Furthermore, a cycle $C$ is {\em prime} if it is not a power of 
a strictly smaller cycle. 
Note that each equivalence class of prime, reduced cycles of a graph $G$ 
corresponds to a unique conjugacy class of 
the fundamental group $ \pi {}_1 (G,u)$ of $G$ at a vertex $u$ of $G$. 
Furthermore, an equivalence class of prime cycles of a graph $G$ 
is called a {\em primitive periodic orbit} of $G$(see \cite{Smi}).   

The {\em Ihara zeta function} of a graph $G$ is 
a function of a complex variable $t$ with $\mid t \mid $ sufficiently small, 
defined by 
\[
{\bf Z} (G, t)= {\bf Z}_G (t)= \prod_{[p]} (1- t^{ \mid p \mid } )^{-1} ,
\]
where $[p]$ runs over all equivalence classes of prime, reduced cycles of $G$(see \cite{Ihara}).

\newtheorem{theorem}{Theorem}
\begin{theorem}[Ihara; Bass]
Let $G$ be a connected graph. 
Then the reciprocal of the Ihara zeta function of $G$ is given by 
\[
{\bf Z} (G, t)^{-1} =(1- t^2 )^{r-1} \det ( {\bf I} -t {\bf A} (G)+ 
t^2 ({\bf D} -{\bf I} )), 
\]
where $r$ and ${\bf A} (G)$ are the Betti number and the adjacency matrix 
of $G$, respectively, and ${\bf D} =( d_{ij} )$ is the diagonal matrix 
with $d_{ii} = v_i = \deg u_i $ where $V(G)= \{ u_1 , \cdots , u_n \} $. 
\end{theorem}

Let $G$ be a connected graph and $V(G)= \{ u_1 , \cdots , u_n \}$. 
Then we consider an $n \times n$ matrix 
${\bf W} =( w_{ij} )_{1 \leq i,j \leq n }$ with $ij$ entry 
the complex variable $w_{ij}$ if $( u_i , u_j ) \in D(G)$, 
and $w_{ij} =0$ otherwise. 
The matrix ${\bf W} = {\bf W} (G)$ is called the 
{\em weighted matrix} of $G$.
Furthermore, let $w( u_i , u_j )= w_{ij}, \  u_i , u_j \in V(G)$ and 
$w(b)= w_{ij}, b=( u_i , u_j ) \in D(G)$. 
For each path $P=( e_{i_1} , \cdots , e_{i_r} )$ of $G$, the {\em norm} 
$w(P)$ of $P$ is defined as follows: 
$w(P)= w( e_{ i_1 }) w(e_{ i_2 } ) \cdots  w( e_{ i_r }) $.

Let $G$ be a connected graph with $n$ vertices and $m$ unoriented edges, 
and ${\bf W} = {\bf W} (G)$ a weighted matrix of $G$.
Two $2m \times 2m$ matrices 
${\bf B} = {\bf B} (G)=( {\bf B}_{e,f} )_{e,f \in R(G)} $ and 
${\bf J}_0 ={\bf J}_0 (G) =( {\bf J}_{e,f} )_{e,f \in R(G)} $ 
are defined as follows: 
\[
{\bf B}_{e,f} =\left\{
\begin{array}{ll}
w(f) & \mbox{if $t(e)=o(f)$, } \\
0 & \mbox{otherwise, }
\end{array}
\right.
\  
{\bf J}_{e,f} =\left\{
\begin{array}{ll}
1 & \mbox{if $f= \hat{e} $, } \\
0 & \mbox{otherwise.}
\end{array}
\right.
\]
Then the {\em second weighted zeta function} of $G$ is defined by 
\[
{\bf Z}_1 (G,w,t)= \det ( {\bf I}_n -t ( {\bf B} - {\bf J}_0 ) )^{-1} . 
\]
If $w(e)=1$ for any $e \in D(G)$, then the zeta function of $G$ is the 
Ihara zeta function of $G$.

\begin{theorem}[Sato]
Let $G$ be a connected graph, and 
let ${\bf W} = {\bf W} (G)$ be a weighted matrix of $G$. 
Then the reciprocal of the second weighted zeta function of $G$ is given by 
\[
{\bf Z}_1 (G,w,t )^{-1} =(1- t^2 )^{m-n} 
\det ({\bf I}_n -t {\bf W} (G)+ t^2 ( \tilde{{\bf D}} - {\bf I}_n )) , 
\]
where $n= \mid V(G) \mid $, $m= \mid E(G) \mid $ and 
$\tilde{{\bf D}} =( d_{ij} )$ is the diagonal matrix 
with $d_{ii} = \sum_{o(b)= u_i } w(e)$, $V(G)= \{ u_1 , \cdots , u_n \} $. 
\end{theorem}

Next, we state the bond scattering matrix of a graph. 
Let $G$ be a connected graph with $n$ vertices and $m$ edges, 
$V(G)= \{ u_1 , \ldots , u_n \} $ and $D(G)= \{ b_1 , \ldots , b_m , 
b_{m+1} , \ldots , b_{2m} \} $ such that $b_{m+j} = b^{-1}_j (1 \leq j \leq m)$. 
The {\em Laplacian (matrix)} ${\bf L} = {\bf L} (G)$ of $G$ is defined by 
\[
{\bf L} = {\bf L} (G)=- {\bf A} (G) + {\bf D} . 
\]
Let $\lambda $ be a eigenvalue of ${\bf L}$ and $\psi =(\psi {}_1 , \ldots , \psi {}_n )$ 
the eigenvector corresponding to $\lambda $. 
For each arc $b=(u_j , u_l)$, one associates a {\em bond wave function} 
\[
\psi {}_b (x)= a_b {\rm e}^{i \pi x/4} + a_{b^{-1}} {\rm e}^{-i \pi x /4} , 
\  x= \pm 1 
\]
under the condition 
\[
\psi {}_b (1)= \psi {}_j , \psi {}_b (-1)= \psi {}_l . 
\]
We consider the following three conditions: 
\begin{enumerate} 
\item {\em uniqueness}: The value of the eigenvector at the vertex $u_j $, $\psi {}_j $, 
computed in the terms of the bond wave functions is the same for all the arcs 
emanating from $u_j $.  
\item {\em $\psi $ is an eigenvector of ${\bf L} $}; 
\item {\em consistency}: The linear relation between the incoming and 
the outgoing coefficients (1) must be satisfied simultaneously at all vertices. 
\end{enumerate} 

By the uniqueness, we have 
\[
a_{b_1} {\rm e}^{i \pi /4} +a_{b^{-1}_1 } {\rm e}^{-i \pi /4} =
a_{b_2} {\rm e}^{i \pi /4} +a_{b^{-1}_2 } {\rm e}^{-i \pi /4} = \cdots 
= a_{b_{d_j}} {\rm e}^{i \pi /4} +a_{b^{-1}_{d_j} } {\rm e}^{-i \pi /4} ,  
\]
where $b_1 , b_2 , \ldots , b_{d_j}$ are arcs emanating from $u_j $, and $d_j = \deg u_j $, 
$i= \sqrt{-1} $. 

By the condition 2, we have 
\[
-\sum^{d_j}_{k=1} ( a_{b_k} {\rm e}^{-i \pi /4} +a_{b^{-1}_k } {\rm e}^{i \pi /4} ) 
=( \lambda -v_j ) \frac{1}{v_j } 
\sum^{d_j}_{k=1} ( a_{b_k} {\rm e}^{i \pi /4} +a_{b^{-1}_k } {\rm e}^{-i \pi /4} ) . 
\]
Thus, for each arc $b$ with $o(b)=u_j $, 
\begin{equation}
a_b = \sum_{t(c)= u_j } \sigma {}^{(u_j)}_{b,c} (\lambda ) a_c , 
\end{equation}
where 
\[
\sigma {}^{(u_j)}_{b,c} (\lambda ) 
=i( \delta {}_{b^{-1}, c} - \frac{2}{d_j} \frac{1}{1-i(1- \lambda / d_j ) } ) ,  
\]
and $\delta {}_{b^{-1},  c} $ is the Kronecker delta. 
The {\em bond scattering matrix} ${\bf U} (\lambda )=( U_{ef} )_{e,f \in D(G)} $ 
of $G$ is defined by 
\[
U_{ef} =\left\{
\begin{array}{ll}
\sigma^{(t(f))}_{e,f} & \mbox{if $t(f)=o(e)$, } \\
0 & \mbox{otherwise}
\end{array}
\right. 
\]
By the consistency, we have 
\[
{\bf U} ( \lambda ) {\bf a} = {\bf a} , 
\]
where ${\bf a} = {}^t (a_{b_1} , a_{b_2} , \ldots , a_{b_{2m}} )$. 
This holds if and only if 
\[
\det ( {\bf I}_{2m} - {\bf U} ( \lambda ))=0 . 
\]

\begin{theorem}[Smilansky]
Let $G$ be a connected graph with $n$ vertices and $m$ edges. 
Then the characteristic polynomial of the bond scattering matrix of $G$ 
is given by 
\[
\det ( {\bf I}_{2m} - {\bf U} ( \lambda )) 
= \frac{ 2^m i^n \det ( \lambda {\bf I}_n + {\bf A} (G) - {\bf D} )}
{\prod^n_{j=1} ( d_j -i d_j + \lambda i)}  
= \prod_{[p]} (1- a_p ( \lambda )  ) , 
\]
where $[p]$ runs over all primitive periodic orbits of $G$, and 
\[
a_p( \lambda  )= \sigma {}^{(t(b_n ))} _{b_1, b_n } \sigma {}^{(t(b_{n-1})) }_{b_n, b_{n-1} } \cdots 
\sigma {}^{(t(b_1)) } _{b_2, b_1 } , 
\  p=( b_1 , b_2 , \ldots , b_n ) 
\]
\end{theorem}

Mizuno and Sato \cite{MS} presented another proof for this Smilansky's formula 
by using the determinant expression of the second weighted zeta function of a graph.

\section{The scattering matrix of a graph with respect to a Hermitian matrix}

Let $G$ be a connected graph with $n$ vertices and $m$ edges, 
$V(G)= \{ 1 , \ldots , n \} $ and $D(G)= \{ e_1 , \ldots , e_m , 
e_{m+1} , \ldots , e_{2m} \} $ such that $e_{m+j} = e^{-1}_j (1 \leq j \leq m)$. 
Furthermore, let an Hermitian matrix ${\bf H} = {\bf H} (G)=( H_{uv} )_{u,v \in V(G)} $ be given as follows: 
\[
H_{uv} =\left\{
\begin{array}{ll}
h_{f} {\rm e}^{2 i \gamma {}_{f}} & \mbox{if $f=(u,v) \in D(G)$, } \\
0 & \mbox{otherwise, }
\end{array}
\right. 
\]
where, for each $f \in D(G)$,  
\[
h_f = h_{f^{-1} } \geq 0 \ and \ \gamma {}_f = - \gamma {}_{f^{-1} } \in [ - \pi /2 , \pi /2] . 
\]
If $H_{uv} = H_f$ is real and negative, then we choose $\gamma {}_f = \pi /2$ if $u \geq v$ and $\gamma {}_f =- \pi /2$ if $u<v$. 
Set 
\[
h(u,v)= h_{uv} =h_f \ and \  \gamma (u,v)= \gamma {}_{uv} = \gamma {}_f \  for \ f=(u,v) \in D(G) . 
\]

Now, let $\lambda $ be an eigenvalue of ${\bf H}$ and $\psi =(\psi {}_1 , \ldots , \psi {}_n )$ 
the eigenvector corresponding to $\lambda $. 
For each arc $b=(u,v)$, one associates a {\em bond wave function} 
\[
\psi {}_b (x)= \frac{ {\rm e}^{i \gamma {}_b } }{\sqrt{h_b}} ( a_{b^{-1}} {\rm e}^{i \pi x/4} + a_{b} {\rm e}^{-i \pi x /4} ) , 
\  x= \pm 1 
\]
under the condition 
\[
\psi {}_b (1)= \psi {}_u , \psi {}_b (-1)= \psi {}_v . 
\]
We consider the following three conditions: 
\begin{enumerate} 
\item {\em uniqueness}: The value of the eigenvector at the vertex $u $, $\psi {}_u $, 
computed in the terms of the bond wave functions is the same for all the arcs 
emanating from $u $.  
\item {\em $\psi $ is an eigenvector of ${\bf H} $}; 
\item {\em consistency}: The linear relation between the incoming and 
the outgoing coefficients (1) must be satisfied simultaneously at all vertices. 
\end{enumerate} 

By the uniqueness 1, we have 
\[
\frac{ {\rm e}^{i \gamma {}_{b_1} } }{\sqrt{h_{b_1}} } ( a_{b^{-1}_1} {\rm e}^{i \pi /4} +a_{b_1 } {\rm e}^{-i \pi /4} )=
\frac{ {\rm e}^{i \gamma {}_{b_2} } }{\sqrt{h_{b_2}} } ( a_{b^{-1}_2} {\rm e}^{i \pi /4} +a_{b_2 } {\rm e}^{-i \pi /4} )= \cdots 
\]
\[
= \frac{ {\rm e}^{i \gamma {}_{b_d} } }{\sqrt{h_{b_d}} } ( a_{b^{-1}_d} {\rm e}^{i \pi /4} +a_{b_d } {\rm e}^{-i \pi /4} )=  \psi {}_u ,  
\]
where $b_1 , b_2 , \ldots , b_{d} $ are arcs emanating from $u$, and $d= \deg u$, 
$i= \sqrt{-1} $. 

By the condition 2, we have 
\[ 
(H_{uu} - \lambda ) \psi {}_u + \sum_{v \in {\cal E}_u } H_{uv} \psi {}_v = 0,  
\]
and so, 
\[
(H_{uu} - \lambda ) \frac{ {\rm e}^{i \gamma {}_{b_1} } }{\sqrt{h_{b_1}} } ( a_{b^{-1}_1} {\rm e}^{i \pi /4} +a_{b_1 } {\rm e}^{-i \pi /4} )
=- \frac{1}{d} \sum^{d}_{k=1} H_{b_j } \frac{ {\rm e}^{i \gamma {}_{b_k}} }{\sqrt{h_{b_k}} } 
( a_{b^{-1}_k} {\rm e}^{i \pi /4} +a_{b_k } {\rm e}^{-i \pi /4} ) , 
\]
where ${\cal E}_u = \{ f \in D(G) \mid o(f)=u \} $. 
Thus, for each arc $b$ with $o(b)=u$, 
\[ 
a^{-1}_b = i a_b -2 \sum^d_{k=1} \frac{ \sqrt{h_b} \sqrt{h_{b_k} } }{H_{uu} -\lambda -i \Gamma {}_u } 
{\rm e}^{i( \gamma {}_{b_k } + \gamma {}_{b^{-1}} )} a_{b_k } , 
\]
where 
\[
\Gamma {}_u = \sum^d_{k=1} h_{b_k} . 
\]

Let $e=b^{-1} $, $f=b_k$ and 
\[
\sigma^{(u)}_{ef} ( \lambda )= i \delta {}_{e^{-1} f} -2 \frac{ \sqrt{h_e} \sqrt{h_f} }{H_{uu} -\lambda -i \Gamma {}_u } 
{\rm e}^{i( \gamma {}_{f} + \gamma {}_{e} )} ,  
\] 
where $\delta {}_{e^{-1} f} $ is the Kronecker delta. 
Then we have 
\begin{equation} 
a_e = \sum_{o(f)=u} \sigma^{(u)}_{ef} ( \lambda ) a_f   
\end{equation}
for each arc $e$ such that $t(e)=u$.  
The {\em bond scattering matrix} ${\bf U} (\lambda )=( U_{ef} )_{e,f \in D(G)} $ 
of $G$ is defined by 
\[
U_{ef} =\left\{
\begin{array}{ll}
\sigma^{(t(e))}_{e,f} & \mbox{if $t(e)=o(f)$, } \\
0 & \mbox{otherwise}
\end{array}
\right. 
\]
By the consistency 3, we have 
\[
{\bf U} ( \lambda ) {\bf a} = {\bf a} , 
\]
where ${\bf a} = {}^t (a_{b_1} , a_{b_2} , \ldots , a_{b_{2m}} )$. 
This holds if and only if 
\[
\det ( {\bf I}_{2m} - {\bf U} ( \lambda ))=0 . 
\]

We present another proof of Theorem 4 by using the technique on the Ihara zeta function, which is different from a proof in \cite{GS}.

\begin{theorem}[Gnutzmann and Smilansky]
Let $G$ be a connected graph with $n$ vertices $1, \ldots , n$ and $m$ edges. 
Then, for the bond scattering matrix of $G$,  
\[
\det ( {\bf I}_{2m} - {\bf U} ( \lambda )) 
= \frac{ (-1)^n 2^m \det ( \lambda {\bf I}_n - {\bf H} )}
{\prod^n_{j=1} ( H_{jj} - \lambda -i \Gamma {}_j )} .  
\]
\end{theorem}

{\bf Proof}.  The argument is an analogue of Watanabe and Fukumizu's method \cite{WF}. 

Let $G$ be a connected graph with $n$ vertices and $m$ edges, 
$V(G)= \{ 1 , \cdots , n \}$ 
and $D(G)= \{ b_1 , \ldots , b_m , b^{-1}_1 , \ldots , b^{-1}_m \}$. 
Set $d_j = \deg j $ and 
\[
x_j = \frac{2}{ H_{jj} - \lambda -i \Gamma {}_j } 
\]
for each $j=1, \ldots ,n$. 
Furthermore, for $e \in D(G)$, let \[
w(e)= \sqrt{h_e } {\rm e}^{i \gamma {}_e} . 
\]
Them we have 
\[
\sigma^{(t(e))}_{ef} ( \lambda )= i \delta {}_{e^{-1} f} - x_{t(e)} w(e) w(f) . 
\]

Now, we consider a $2m \times 2m$ matrix 
${\bf B} =( B_{ef} )_{e,f \in D(G)}$ given by 
\[
B_{ef} =\left\{
\begin{array}{ll}
x_{o(f)} w(e) w(f) & \mbox{if $t(e)=o(f)$, } \\
0 & \mbox{otherwise}
\end{array}
\right.
\] 
Let ${\bf K} =( {\bf K}_{i,j} )$ ${}_{1 \leq i \leq 2m; 
1 \leq j \leq n} $ be the $2m \times n$ matrix defined 
as follows: 
\[
{\bf K}_{i,j} :=\left\{
\begin{array}{ll}
x_j w(b_i) & \mbox{if $o( b_i )=j$, } \\
0 & \mbox{otherwise. } 
\end{array}
\right.
\] 
Furthermore, we define two $2m \times n$ matrices 
${\bf L} =( {\bf L}_{i,j} )_{1 \leq i \leq 2m; 
1 \leq j \leq n} $ and 
${\bf M} =( {\bf M}_{i,j} )_{1 \leq i \leq 2m; 
1 \leq j \leq n} $ as follows: 
\[
{\bf L}_{i,j} :=\left\{
\begin{array}{ll}
w(b_i) & \mbox{if $t( b_i )=j$, } \\
0 & \mbox{otherwise, } 
\end{array}
\right.
\ 
{\bf M}_{i,j} :=\left\{
\begin{array}{ll}
w(b_i) & \mbox{if $o( b_i )=j$, } \\
0 & \mbox{otherwise. } 
\end{array}
\right.
\]
Note that 
\begin{equation}
{\bf K} ={\bf M} 
\left[
\begin{array}{ccc}
x_{1} & \  & 0 \\
\  & \ddots & \  \\
0 & \  & x_{n } 
\end{array}
\right]
= {\bf M} {\bf X} . 
\end {equation}

Furthermore, we have 
\begin{equation}
{\bf L} {}^t {\bf K} = {\bf B}  
\end{equation}
and 
\begin{equation}
{}^t {\bf M} {\bf L} = {\bf H} .  
\end{equation}
Note that 
\[
H_{uv} = w(u,v)^2 \ if \ (u,v) \in D(G) . 
\]

But, since 
\[
U_{ef} =\left\{
\begin{array}{ll}
-x_{t(e)} w(e)w(f) & \mbox{if $t(e)=o(f)$ and $f \neq e^{-1} $, } \\
i-x_{t(e)} w(e)w(f) & \mbox{if $f= e^{-1} $, } \\
0 & \mbox{otherwise, } 
\end{array}
\right.
\]
we have 
\[
{\bf U} ( \lambda ) = i {\bf J}_0 - {\bf B} . 
\]
Furthermore, if ${\bf A}$ and ${\bf B}$ are an $r \times s$ and an $s \times r$ 
matrix, respectively, then we have 
\[
\det ( {\bf I}_{r} - {\bf A} {\bf B} )= 
\det ( {\bf I}_s - {\bf B} {\bf A} ) . 
\] 
Thus, 
\[
\begin{array}{rcl} 
\det ( {\bf I}_{2m} -u {\bf U} ( \lambda ) ) & = & \det ( {\bf I}_{2m} -u (i {\bf J}_0 - {\bf B} )) \\ 
\  &   &                \\ 
\  & = & \det ({\bf I}_{2m} -iu {\bf J}_0 +u {\bf L} \  {}^t {\bf K} ) \\ 
\  &   &                \\ 
\  & = & \det ({\bf I}_{2m} +u {\bf L} \  {}^t {\bf K} ( {\bf I}_{2m} -iu {\bf J}_0)^{-1} ) \det ( {\bf I}_{2m} -iu {\bf J}_0 ) \\ 
\  &   &                \\ 
\  & = & \det ({\bf I}_{n} +u \  {}^t {\bf K} ( {\bf I}_{2m} -iu {\bf J}_0)^{-1} {\bf L} ) \det ( {\bf I}_{2m} -iu {\bf J}_0 ) . 
\end{array}
\]

Arrange arcs of $D(G)$ as follows: $b_1 , b^{-1}_1 , \ldots , b_m , b^{-1}_m $. 
Then we have 
\[ 
\det ( {\bf I}_{2m} -iu {\bf J}_0 )
=\det (  \left[
\begin{array}{cccc}
1 & -iu & \ldots & 0 \\
-iu & 1 &  &  \\ 
\vdots &  & \ddots & \\
0 &  &  &   
\end{array}
\right] 
)=(1+ u^2 )^m . 
\]
Furthermore, 
\[
\begin{array}{rcl} 
( {\bf I}_{2m} -iu {\bf J}_0 )^{-1} & = & \left[
\begin{array}{cccc}
1 & -iu & \ldots & 0 \\
-iu & 1 &  &  \\ 
\vdots &  & \ddots & \\
0 &  &  &   
\end{array}
\right]^{-1} 
\\
\  &   &                \\ 
\  & = & \frac{1}{1+ u^2 } 
\left[
\begin{array}{cccc}
1 & iu & \ldots & 0 \\
iu & 1 &  &  \\ 
\vdots &  & \ddots & \\
0 &  &  &   
\end{array}
\right] 
\\ 
\  &   &                \\ 
\  & = & \frac{1}{1+ u^2 } ( {\bf I}_{2m} +iu {\bf J}_0 ) . 
\end{array}
\]
Therefore, it follows that 
\[
\begin{array}{rcl} 
\  &   & \det ( {\bf I}_{2m} -u {\bf U} ( \lambda ) ) \\ 
\  &   &                \\ 
\  & = & \det ({\bf I}_{n} + \frac{u}{1+ u^2 } \  {}^t {\bf K} ( {\bf I}_{2m} +iu {\bf J}_0) {\bf L} )(1+ u^2 )^m  \\ 
\  &   &                \\ 
\  & = & (1+ u^2 )^{m-n} \det ((1+ u^2 ) {\bf I}_{n} +u  \ {}^t {\bf K} {\bf L} +i u^2 \ {}^t {\bf K} {\bf J}_0 {\bf L} ) . 
\end{array}
\]

But, we have 
\[
{}^t {\bf K} {\bf L} ={\bf X} {}^t {\bf M} {\bf L} = {\bf X} {\bf H} . 
\]
Furthermore, we have 
\[
{}^t {\bf K} {\bf J}_0 {\bf L} ={\bf X} {}^t {\bf M} {\bf J}_0 {\bf L} . 
\]
Then, for $u,v \in V(G)$, we have  
\[
\begin{array}{rcl} 
\  &  & ( {}^t {\bf M} {\bf J}_0 {\bf L} )_{uv} \\
\  &   &                \\ 
\  & = & \delta {}_{uv} \sum_{o(e)=u} ( {}^t {\bf M} )_{ue} ( {\bf J}_0 )_{e e^{-1} } ({\bf L} )_{e^{-1} v} \\ 
\  &   &                \\ 
\  & = & \delta {}_{uv} \sum_{o(e)=u} w(e) \cdot 1 \cdot w(e^{-1} ) \\
\  &   &                \\ 
\  & = & \delta {}_{uv} \sum_{o(e)=u} \sqrt{h_e} {\rm e}^{i \gamma {}_e } \sqrt{h_e} {\rm e}^{-i \gamma {}_e } \\ 
\  &   &                \\ 
\  & = & \delta {}_{uv} \sum_{o(e)=u} h_e =  \delta {}_{uv} \Gamma {}_u . 
\end{array}
\]

Now, let 
\[
{\bf D}_L = \left[
\begin{array}{ccc}
\Gamma {}_{1} & \  & 0 \\
\  & \ddots & \  \\
0 & \  & \Gamma {}_{n } 
\end{array}
\right]
. 
\]
Then 
\[
{}^t {\bf K} {\bf J}_0 {\bf L} ={\bf X} {\bf D}_L . 
\]
Thus, 
\[
\det ( {\bf I}_{2m} -u {\bf U} ( \lambda ) )=(1+ u^2 )^{m-n} \det ((1+ u^2 ) {\bf I}_{n} +u {\bf X} {\bf H} +i u^2 {\bf X} {\bf D}_L ) . 
\]
Substituting $u=1$, we obtain 
\[ 
\begin{array}{rcl} 
\  &  & \det ( {\bf I}_{2m} - {\bf U} ( \lambda )) \\ 
\  &   &                \\ 
\  & = & 2^{m-n} \det (2 {\bf I}_n + {\bf X} {\bf H}+i {\bf X} {\bf D}_L ) \\ 
\  &   &                \\ 
\  & = & 2^{m-n} \det ( 
\left[
\begin{array}{ccccc}
\ldots & 2+ i \Gamma {}_u \frac{2}{H_{uu} - \lambda -i \Gamma {}_u } & \ldots 
& \frac{2}{H_{uu} - \lambda -i \Gamma {}_u } h_{uv} {\rm e}^{2i \gamma {}_{uv} } & \ldots 
\end{array}
\right] 
) 
\\ 
\  &   &                \\ 
\  & = & \frac{ 2^{m}}{ \prod^n_{u=1} ( H_{uu} - \lambda -i \Gamma {}_u )} \det(- \lambda {\bf I}_n + {\bf H} ) \\
\  &   &                \\ 
\  & = & \frac{ (-1)^n 2^{m}}{ \prod^n_{u=1} ( H_{uu} - \lambda -i \Gamma {}_u )} \det( \lambda {\bf I}_n - {\bf H} ) . 
\end{array}
\]
$\Box$

\section{The Euler product with respect to the scattering matrix}

We present the Euler product for the determinant formula of the scattering matrix ${\bf U} ( \lambda )$ of a graph.

\begin{theorem}
Let $G$ be a connected graph with $m$ edges, and ${\bf H} = {\bf H} (G)=( H_{uv} )_{u,v \in V(G)} $ 
an Hermitian matrix defined in Section 2.  
Then the characteristic polynomial of the bond scattering matrix of $G$ induced from ${\bf H} $ 
is given by 
\[
\det ( {\bf I}_{2m} -u {\bf U} ( \lambda )) = \prod_{[C]} (1- w_C u^{|C|} ) , 
\]
let ${c} $ runs over all equivalence classes of prime cycles in $G$, and 
\[
w_C= \sigma {}^{(t(e_1 ))}_{e_1 e_2 } \sigma {}^{(t(e_2 ))}_{e_2 e_3 } \cdots \sigma {}^{(t(e_n ))}_{e_n e_1 } , 
\  C=( b_1 , b_2 , \ldots , b_n ) 
\]
\end{theorem}

{\bf Proof}.  Let $D(G)= \{ b_1 , \cdots , b_{2m} \} $ such that 
$b_{m+j} = b^{-1}_j (1 \leq j \leq m)$. 
Set $ {\bf U} = {\bf U} ( \lambda ) $. 
Since 
\[
\log \det ( {\bf I} - u {\bf F} )= {\rm Tr} \log ( {\bf I} - u {\bf F} ), 
\]
for a square matrix ${\bf F} $, we have 
\[
\log \det ( {\bf I} - u {\bf U} )= {\rm Tr} \log ( {\bf I} - u {\bf U} )
= - \sum^{\infty}_{ k=1} \frac{ {\rm Tr} ( {\bf U}^k )}{k} u^k . 
\]
Here, 
\[
{\rm Tr} ( {\bf U}^k )= \sum_C w_C , 
\]
where $C$ runs over all cycles of length $k$ in $G$, and 
\[
w_C= \sigma {}^{(t(e_1 ))}_{e_1 e_2 } \sigma {}^{(t(e_2 ))}_{e_2 e_3 } \cdots \sigma {}^{(t(e_k ))}_{e_k e_1 } , 
\  C=( b_1 , b_2 , \ldots , b_k ) 
\]
Thus, 
\[
\begin{array}{rcl} 
u \frac{d}{du} \log \det ( {\bf I}_{2m} -u {\bf U} ) & = & \sum^{\infty}_{k=1} {\rm Tr} ( {\bf U}^k ) u^k \\ 
\  &   &                \\ 
\  & = & \sum_C w_C u^{|C|} , 
\end{array}
\]
where $C$ runs over all cycles in $G$.

Now, let $C$ be any cycle in $G$. 
Then there exists exactly one prime cycle $D$ such that 
\[
C=D^l . 
\]
Thus, we have 
\[
u \frac{d}{du} \log \det ( {\bf I}_{2m} -u {\bf U} )= - \sum_D \sum^{\infty}_{k=1} w^k_D u^{k |D|} , 
\]
and so, 
\[
\frac{d}{du} \log \det ( {\bf I}_{2m} -u {\bf U} )= - \sum_D \sum^{\infty}_{k=1} w^k_D u^{k |D|-1} ,  
\]
where $D$ runs over all prime cycles in $G$. 
Therefore, it follows that 
\[
\begin{array}{rcl} 
\log \det ( {\bf I}_{2m} -u {\bf U} ) & = & - \sum_D \sum^{\infty}_{k=1} \frac{w^k_D}{k|D|} u^{k |D|} \\ 
\  &   &                \\ 
\  & = & - \sum_{[D]} \sum^{\infty}_{k=1} \frac{|D|}{k|D|} w^k_D u^{k |D|} \\ 
\  &   &                \\ 
\  & = & - \sum_{[D]} \sum^{\infty}_{k=1} \frac{1}{k} w^k_D u^{k |D|} \\
\  &   &                \\ 
\  & = & \sum_{[D]} \log (1- w_D u^{|D|} ) . 
\end{array}
\]
Hence, 
\[
\det ( {\bf I}_{2m} -u {\bf U} ( \lambda )) = \prod_{[C]} (1- w_C u^{|C|} ) , 
\]
$\Box$

\section{Scattering matrix of a regular covering of a graph}

Let $G$ be a connected graph, and 
let $N(v)= \{ w \in V(G) \mid (v,w) \in D(G) \} $ denote the 
neighbourhood  of a vertex $v$ in $G$. 
A graph $H$ is a {\em covering} of $G$ 
with projection $ \pi : H \longrightarrow G $ if there is a surjection
$ \pi : V(H) \longrightarrow V(G)$ such that
$ \pi {\mid}_{N(v')} : N(v') \longrightarrow N(v)$ is a bijection 
for all vertices $v \in V(G)$ and $v' \in {\pi}^{-1} (v) $.
When a finite group $\Pi$ acts on a graph $G$, 
the {\em quotient graph} $G/ \Pi$ is a graph 
whose vertices are the $\Pi$-orbits on $V(G)$, 
with two vertices adjacent in $G/ \Pi$ if and only if some two 
of their representatives are adjacent in $G$.
A covering $ \pi : H \longrightarrow G$ is {\em regular} 
if there is a subgroup {\it B} of the automorphism group $Aut \  H$ 
of $H$ acting freely on $H$ such that the quotient graph $H/ {\it B} $ 
is isomorphic to $G$.

Let $G$ be a graph and $ \Gamma $ a finite group.
Then a mapping $ \alpha : D(G) \longrightarrow \Gamma $
is an {\em ordinary voltage} {\em assignment}
if $ \alpha (v,u)= \alpha (u,v)^{-1} $ for each $(u,v) \in D(G)$.
The pair $(G, \alpha )$ is an {\em ordinary voltage graph}.
The {\em derived graph} $G^{ \alpha } $ of the ordinary
voltage graph $(G, \alpha )$ is defined as follows:
$V(G^{ \alpha } )=V(G) \times \Gamma $ and $((u,h),(v,k)) \in 
D(G^{ \alpha })$ if and only if $(u,v) \in D(G)$ and $k=h \alpha (u,v) $. 
The {\em natural projection} 
$ \pi : G^{ \alpha } \longrightarrow G$ is defined by 
$ \pi (u,h)=u$. 
The graph $G^{ \alpha }$ is a {\em derived graph covering} of $G$ 
with voltages in $ \Gamma $ or a {\em $ \Gamma $-covering} of $G$. 
Note that $| {\cal E}_{(u,h)} |=| {\cal E}_u |$ for each $(u,h) \in V(G^{\alpha } )$. 
The natural projection $ \pi $ commutes with the right 
multiplication action of the $ \alpha (e), e \in D(G)$ and 
the left action of $ \Gamma $ on the fibers: 
$g(u,h)=(u,gh), g \in \Gamma $, which is free and transitive. 
Thus, the $ \Gamma $-covering $G^{ \alpha }$ is a $ \mid \Gamma \mid $-fold
regular covering of $G$ with covering transformation group $ \Gamma $.
Furthermore, every regular covering of a graph $G$ is a 
$ \Gamma $-covering of $G$ for some group $ \Gamma $ (see \cite{GT}).

Let $G$ be a connected graph, $ \Gamma $ be a finite group and 
$ \alpha : D(G) \longrightarrow \Gamma $ be an ordinary voltage assignment. 
In the $\Gamma $-covering $G^{ \alpha } $, set $v_g =(v,g)$ and $e_g =(e,g)$, 
where $v \in V(G), e \in D(G), g\in \Gamma $. 
For $e=(u,v) \in D(G)$, the arc $e_g$ emanates from $u_g$ and 
terminates at $v_{g \alpha (e)}$. 
Note that $ e^{-1}_g =(e^{-1} )_{g \alpha (e)}$. 

Let $G$ be a connected graph, $ \Gamma $ be a finite group and 
$ \alpha : D(G) \longrightarrow \Gamma $ be an ordinary voltage assignment. 
Furthermore, let ${\bf H} = {\bf H} (G)=( H_{uv} )_{u,v \in V(G)} $ be an Hermitian matrix such that 
\[
H_{uv} =\left\{
\begin{array}{ll}
h_{f} {\rm e}^{2 i \gamma {}_{f}} & \mbox{if $f=(u,v) \in D(G)$, } \\
0 & \mbox{otherwise, }
\end{array}
\right. 
\]
where, for each $f \in D(G)$,  
\[
h_f = h_{f^{-1} } \geq 0 \ and \ \gamma {}_f = - \gamma {}_{f^{-1} } \in [ - \pi /2 , \pi /2] . 
\] 
We give the function $\tilde{h} : D(G^{\alpha } ) \longrightarrow \mathbb{R} $ and 
$\tilde{ \gamma } : D(G^{\alpha } ) \longrightarrow [ -\pi /2, \pi /2] $ induced from $h$ and $\gamma $, respectively, 
as follows: 
\[
\tilde{h} (u_g , v_k )= h_{uv} \ and \ \tilde{ \gamma } (u_g , v_k )= \gamma {}_{uv} 
\ if \ (u,v) \in D(G) \ and \ k=g \alpha (u,v) . 
\]
Furthermore, we consider the Hermitian matrix $\tilde{{\bf H}} ={\bf H} ( G^{\alpha } )=
( H_{u_g v_k } )_{ u_g v_k \in V(G^{\alpha } )} $ of $G^{\alpha } $ induced from ${\bf H} $. 
At first, let 
\[
H_{u_g u_g } = H_{uu} \ for \ each \ g \in \Gamma . 
\]
For $(u_g , v_k ) \in D(G^{\alpha } )$, we have 
\[
H_{u_g v_k } = \tilde{h}(u_g , v_k ) e^{2i \tilde{\gamma} (u_g , v_k )} = h_{uv} e^{2i \gamma {}_{uv}} . 
\]
Thus, 
\[ 
H_{u_g v_k } =\left\{
\begin{array}{ll}
h_{uv} e^{2i \gamma {}_{u_g v_k}} & \mbox{if $(u,v) \in D(G)$ and $k=g \alpha (u,v)$, } \\
H_{uu} & \mbox{if $u=v$ and $k=g$, } \\ 
0 & \mbox{otherwise. }
\end{array}
\right.
\]

Next, we consider the bond wave function of the regular covering $G^{\alpha }$ of $G$. 
Let $V(G)= \{ v_1 , \ldots , v_n \} $, 
$D(G)= \{ e_1 , \ldots , e_m , e^{-1}_1 , \ldots , e^{-1}_m \} $ and 
$\Gamma = \{ g_1 =1 , g_2 , \ldots , g_p \} $. 
Let $\lambda $ be a eigenvalue of 
$\tilde{{\bf H}} = {\bf H} (G^{\alpha })$, and let 
$\tilde{\phi} =(\phi {}_{v_1 , g_1 }, \ldots , \phi {}_{v_1 , g_p }, \ldots, 
\phi {}_{v_n , g_1 }, \ldots , \phi {}_{v_n , g_p } )$ be the eigenvector corresponding 
to $\lambda $, where  $\phi {}_{v_i , g_j } $ corresponds to the vertex 
$(v_i ,g_j ) \  (1 \leq i \leq n; 1 \leq j \leq p)$ of $G^{\alpha } $.   
Furthermore let $b_g =(v_g , z_{g \alpha (b)} )$ be any arc of $G^{\alpha } $, where 
$b=(v,z) \in D(G)$, $g \in \Gamma $. 
Then the bond wave function of $G^{\alpha } $ is  
\[
\phi {}_{b_g} (x)= \frac{e^{i \gamma {}_b }}{ \sqrt{h_b}} ( a_{b^{-1}_g} {\rm e}^{i \pi x/4} + a_{b_g} {\rm e}^{-i \pi x /4} ) ,   
\  x= \pm 1 , \ i= \sqrt{-1}  
\]
under the condition 
\[
\phi {}_{b_g} (1)= \phi {}_{v_g} \  and \   \phi {}_{b_g} (-1)= \phi {}_{z_{g \alpha (b)}} . 
\]

By (1), we have  
\[
\begin{array}{rcl} 
a_{b^{-1}_g} & = & i \delta {}_{b^{-1}_g e_g} -2 \sum_{o(e_g )= v_g } 
\frac{ \sqrt{\tilde{h}_{e_g }} \sqrt{\tilde{h}_{b_g }} }{H_{v_g v_g } - \lambda -i \Gamma {}_{v_g } } 
e^{i ( \tilde{\gamma }_{b_g } + \tilde{\gamma }_{e_g })} a_{e_g } \\
\  &   &                \\ 
\  & = & \sum_{o(e_g )= v_g } \sigma^{(v_g )}_{b_g e_g } a_{e_g } 
\end{array}
\] 
for each arc $b_g$ with $o(b_g)=v_g $, 
where 
\[
\sigma {}^{(v_g)}_{b_g e_g} =i \delta {}_{b^{-1}_g e_g} 
-2 \frac{\sqrt{\tilde{h}_{e_g }} \sqrt{\tilde{h}_{b_g }} }{H_{v_g v_g } - \lambda -i \Gamma {}_{v_g } } 
e^{i ( \tilde{\gamma }_{b_g } + \tilde{\gamma }_{e_g })} a_{e_g }  
\] 
and 
\[
\tilde{h}_{e_g} = \tilde{h} (e_g ), \ \tilde{ \gamma }_{e_g} = \tilde{ \gamma} (e_g ) . 
\]
By the definitions of $\tilde{h}$, $ \tilde{\gamma } $ and $ \tilde{{\bf H}} $, we have 
\[
\sigma {}^{(v_g)}_{b_g e_g} =i \delta {}_{b^{-1} e} 
-2 \frac{\sqrt{h_e } \sqrt{h_b} }{H_{vv} - \lambda -i \Gamma {}_{v} } 
e^{i ( \gamma {}_{b} + \gamma {}_e)} = \sigma^{(v)}_{be} =\sigma^{(t(b))}_{be} . 
\]
Note that ${\cal E}_{(v,g)} = {\cal E}_v $.  
Thus,  
\[
a_{b^{-1}_g} = \sum_{o(e_g )= v_g } \sigma^{(t(b))}_{be} a_{e_g } . 
\] 
Therefore, the bond scattering matrix 
$\tilde{{\bf U}} ( \lambda )=( U(e_g , f_h) )_{e_g ,f_h \in D(G^{\alpha } )} $ 
of $ G^{ \alpha } $ is given by 
\[
U(e_g , f_h) = 
\left\{
\begin{array}{ll}
\sigma^{(t(e))}_{ef} & \mbox{if $t(f_h)=o(e_g)$, } \\
0 & \mbox{otherwise. }
\end{array}
\right.
\]

But, we have  
\[
x_{v_g } =\frac{2}{ H_{vv} - \lambda -i \Gamma {}_v } = x_v 
\]
for $v_g \in V(G^{\alpha } )$. 
Furthermore, let $\tilde{w} : D(G^{\alpha } ) \longrightarrow {\bf C} $ be given as follows: 
\[
\tilde{w} (e_g )= \sqrt{\tilde{h}_{e_g} } e^{i \tilde{\gamma}_{e_g} } \ for \ each \ e_g \in D(G^{\alpha } ) . 
\]
Then we have 
\[
\tilde{w} (e_g )= \sqrt{h_e } e^{i \gamma {}_e} = w(e), \ e_g \in D(G^{\alpha} ) . 
\]

For $ g \in \Gamma $, let the matrix $ {\bf H} {}_g =( H^{(g)}_{uv} )$ 
be defined by
\[
H^{(g)}_{uv} =\left\{
\begin{array}{ll}
h_{uv} e^{2i \gamma {}_{uv}} & \mbox{if $ \alpha (u,v)=g$ and $(u,v) \in D(G)$, } \\
0 & \mbox{otherwise.}
\end{array}
\right.
\] 
Furthermore, let ${\bf U}_g =( U^{(g)} (e,f) )$ be given by 
\[
U^{(g)} (e,f) =\left\{
\begin{array}{ll}
\sigma^{(t(e))}_{ef} & \mbox{if $t(e)=o(f)$ and $\alpha (e)=g$, } \\
0 & \mbox{otherwise,}
\end{array}
\right.
\] 
Let ${\bf M}_{1} \oplus \cdots \oplus {\bf M}_{s} $ be the block diagonal sum 
of square matrices ${\bf M}_{1} , \ldots , {\bf M}_{s} $. 
If \( {\bf M}_{1} = {\bf M}_{2} = \cdots = {\bf M}_{s} = {\bf M} \),
then we write 
\( s \circ {\bf M} = {\bf M}_{1} \oplus \cdots \oplus {\bf M}_{s} \).
The {\em Kronecker product} $ {\bf A} \bigotimes {\bf B} $
of matrices {\bf A} and {\bf B} is considered as the matrix 
{\bf A} having the element $a_{ij}$ replaced by the matrix $a_{ij} {\bf B}$.

\begin{theorem} 
Let $G$ be a connected graph with $n$ vertices $v_1 , \ldots v_n $ and $m$ 
unoriented edges, $ \Gamma $ be a finite group and 
$ \alpha : D(G) \longrightarrow \Gamma  $ be an ordinary voltage assignment. 
Set $ \mid \Gamma \mid =p$. 
Furthermore, let $ {\rho}_{1} =1, {\rho}_{2} , \cdots , {\rho}_{k} $
be the irreducible representations of $ \Gamma $, and 
$f_i$ be the degree of $ {\rho}_{i} $ for each $i$, where $f_1=1$.

If the $ \Gamma $-covering $G^{ \alpha } $ of $G$ is connected, then, 
for the bond scattering matrix of $G {}^{ \alpha } $,   
\[
\det ( {\bf I}_{2 mp} - \tilde{{\bf U}} ( \lambda ))= \det( {\bf I}_{2m} - {\bf U} ( \lambda )) 
\prod^{k}_{i=2} \det ( {\bf I}_{2m f_i } - \sum_{h} {\rho}_i (h) \bigotimes {\bf U} {}_h )^{f_i} 
\]
\[
= \frac{ 2^{mp} (-1)^{np} \det ( \lambda {\bf I}_{n} - {\bf H} )}{\prod_{u \in V(G)} ( H_{uu} - \lambda -i \Gamma {}_u } 
\prod^{k}_{i=2} \det ( \lambda {\bf I}_{n f_i } - \sum_{h \in \Gamma } {\rho}_{i} (h) \bigotimes {\bf H} {}_{h} 
- {\bf I}_{f_i} \bigotimes {\rm diag} ( {\bf H} ))^{f_i} ,  
\]
where 
\[
{\rm diag} ( {\bf H} ) = \left[
\begin{array}{ccc}
H_{v_1 v_1 } &  & {\bf 0} \\
  & \ddots &  \\
{\bf 0} &  & H_{v_n v_n } 
\end{array}
\right] 
.
\]
\end{theorem}

{\bf Proof }. Let $ \mid \Gamma \mid =p$. 
By Theorem 4, for the bond scattering matrix of $G^{\alpha } $, we have  
\[
\det ( {\bf I}_{2 mp} - \tilde{{\bf U}} ( \lambda )) 
= \frac{ 2^{mp} (-1)^{np} \det ( \lambda {\bf I}_{np} - {\bf H} (G^{\alpha } ))}{\prod_{u \in V(G)} ( H_{uu} - \lambda -i \Gamma {}_u )^p } . 
\]

Let $D(G)= \{ e_1, \ldots , e_m , e_{m+1}, \ldots , e_{2m} \} $ such that $e_{m+j} = e^{-1}_j (1 \leq j \leq m)$ and 
$ \Gamma = \{ 1=g_1, g_2, \ldots ,g_p \} $.
Arrange arcs of $G^{ \alpha } $ in $p$ blocks:
$(e_1,1), \ldots , (e_{2m},1);(e_1,g_2), \ldots , (e_{2m},g_2); \ldots ;$ 
$ (e_1,g_p), \ldots ,(e_{2m},g_p). $
We consider the matrix $\tilde{{\bf U}} ( \lambda ) $
under this order.
For $h \in \Gamma $, the matrix ${\bf P}_{h}=(p^{(h)}_{ij} )$ 
is defined as follows:
\[
p^{(h)}_{ij} = \left\{
\begin{array}{ll}
1 & \mbox{if $g_i h=g_j$,} \\
0 & \mbox{otherwise.}
\end{array}
\right.
\] 

Suppose that $p^{(h)}_{ij} =1 $, i.e., $g_j=g_ih$.
Then $U(e_{g_i}, f_{g_j } ) \neq 0$ if and only if $t(e,g_j)=o(f,g_i)$. 
Furthermore, $t(e,g_j)=o(f,g_i)$ if and only if 
$(o(f), g_{j} )=o(f, g_{j} )=t(e, g_{i} )=(t(e), g_{i} \alpha (e))$. 
Thus, $t(e)=o(f)$ and $ \alpha (e)=g^{-1}_{i} g_j =g^{-1}_{i} g_i h=h$.
Thus, we have   
\[
\tilde{{\bf U}} ( \lambda )= \sum_{h \in \Gamma } {\bf P}_{h} \bigotimes {\bf U} {}_{h} .
\] 
Furthermore, we have 
\[
{\rm diag} ( {\bf H} (G^{\alpha } ))= {\bf I}_p \bigotimes {\rm diag} ( {\bf H} ) .  
\]

Let $\rho$ be the right regular representation of $ \Gamma $.
Furthermore, let $ {\rho}_{1} =1, {\rho}_{2} , \ldots , {\rho}_{k} $
be all inequivalent irreducible representations of $ \Gamma $, and 
$f_i$ the degree of $ {\rho}_{i} $ for each $i$, where 
$f_1=1$.
Then we have $\rho (g)= {\bf P}_{g} $ for $g \in \Gamma $.
Furthermore, there exists a nonsingular matrix ${\bf P}$ such that
${\bf P}^{-1} \rho (g) {\bf P} = (1) \oplus f_2 \circ {\rho}_{2} (g) 
\oplus \cdots \oplus f_k \circ {\rho}_{k} (g)$ 
for each $g \in \Gamma $(see [12]). 
Thus, we have 
\[
{\bf P}^{-1} {\bf P}_g {\bf P} = (1) \oplus f_2 \circ {\rho}_{2} (g) 
\oplus \cdots \oplus f_k \circ {\rho}_{k} (g) . 
\]

Putting 
${\bf F} =( {\bf P}^{-1} \bigotimes {\bf I}_{2m} ) 
\tilde{{\bf U}} ( \lambda ) ( {\bf P} \bigotimes {\bf I}_{2m} )$,
we have 
\[
{\bf F}= \sum_{g \in \Gamma } 
\{ (1) \oplus f_2 \circ {\rho}_{2} (g) \oplus \cdots \oplus 
f_k \circ {\rho}_{k} (g) \} \bigotimes {\bf U} {}_g .
\]
Note that ${\bf U} ( \lambda )= \sum_{g \in \Gamma } {\bf U} {}_g $ and 
$1+ f^2_2 + \cdots + f^2_k =p$.
Therefore it follows that 
\[
\det ( {\bf I}_{2 mp} - \tilde{{\bf U}} ( \lambda ))= 
\det( {\bf I}_{2m} - {\bf U} ( \lambda )) 
\prod^{k}_{i=2} \det ( {\bf I}_{2m f_i } 
- \sum_{g} {\rho}_i (g) \bigotimes {\bf U} {}_g )^{f_i} . 
\]

Next, let $V(G)= \{ v_1, \ldots , v_n \} $. 
Arrange vertices of $G^{ \alpha } $ in $p$ blocks:
$(v_1,1), \ldots , (v_n,1);$ 
$(v_1,g_2), \ldots , (v_n,g_2); 
\ldots ; (v_1,g_p), \ldots ,(v_n,g_p). $
We consider the matrix ${\bf H} ( G {}^{ \alpha } ) $
under this order.

Suppose that $p^{(h)}_{ij} =1 $, i.e., $g_j=g_ih$.
Then $((u,g_i),(v,g_j)) \in D(G {}^{ \alpha } ) $
if and only if $(u,v) \in D(G)$ and 
$g_{j} = g_{i} \alpha (u,v)$.  
If $g_{j} = g_{i} \alpha (u,v)$, then 
$ \alpha (u,v)=g^{-1}_{i} g_j =g^{-1}_{i} g_i h=h$.
Thus we have
\[
{\bf H} (G {}^{ \alpha } )= \sum_{h \in \Gamma } {\bf P}_{h} 
\bigotimes {\bf H} {}_h + {\bf I}_p \bigotimes {\rm diag} ( {\bf H} ) .
\]

Putting 
${\bf E} =( {\bf P}^{-1} \bigotimes {\bf I}_n ) {\bf H} 
(G {}^{ \alpha } ) ( {\bf P} \bigotimes {\bf I}_n )$,
we have 
\[
{\bf E}= \sum_{h \in \Gamma } 
\{ (1) \oplus f_2 \circ {\rho}_{2} (h) \oplus \cdots \oplus 
f_k \circ {\rho}_{k} (h) \} \bigotimes {\bf H} {}_h + {\bf I}_p \bigotimes {\rm diag} ( {\bf H} ) .
\]
Note that ${\bf H} (G) = \sum_{h \in \Gamma } {\bf H} {}_h + {\rm diag} ( {\bf H} )$.  
Therefore it follows that 
\[
\begin{array}{rcl} 
\det ( \lambda {\bf I}_{np} - {\bf H} (G^{\alpha } ))
\  & = & \det (\lambda {\bf I}_{n} - {\bf H} (G)) \\
\  &   &                \\ 
\  & \times & \prod^{k}_{i=2} 
\det ( \lambda {\bf I}_{n f_i } - \sum_{h \in \Gamma } {\rho}_{i} (h) \bigotimes {\bf H} {}_{h} 
- {\bf I}_{f_i} \bigotimes {\rm diag} {\bf H} )^{f_i} .
\end{array}
\] 
Hence, 
\[
\det ( {\bf I}_{2 mp} - \tilde{{\bf U}} ( \lambda ))= \det( {\bf I}_{2m} - {\bf U} ( \lambda )) 
\prod^{k}_{i=2} \det ( {\bf I}_{2m f_i } - \sum_{h} {\rho}_i (h) \bigotimes {\bf U} {}_h )^{f_i} 
\]
\[
= \frac{ 2^{mp} (-1)^{np} \det ( \lambda {\bf I}_{n} - {\bf H} (G))}{\prod_{u \in V(G)} ( H_{uu} - \lambda -i \Gamma {}_u )^p } 
\prod^{k}_{i=2} \det ( \lambda {\bf I}_{n f_i } - \sum_{h \in \Gamma } {\rho}_{i} (h) \bigotimes {\bf H} {}_{h} 
- {\bf I}_{f_i} \bigotimes {\rm diag} {\bf H} )^{f_i} .  
\]
$\Box$

\section{$L$-functions of graphs}

Let $G$ be a connected graph with $n$ vertices and $m$ unoriented edges, 
$ \Gamma $ be a finite group and 
$ \alpha : D(G) \longrightarrow \Gamma $ be an ordinary voltage assignment. 
Furthermore, let ${\bf H} = {\bf H} (G)=( H_{uv} )_{u,v \in V(G)} $ be an Hermitian matrix such that 
\[
H_{uv} =\left\{
\begin{array}{ll}
h_{f} {\rm e}^{2 i \gamma {}_{f}} & \mbox{if $f=(u,v) \in D(G)$, } \\
0 & \mbox{otherwise, }
\end{array}
\right. 
\]
where, for each $f \in D(G)$,  
\[
h_f = h_{f^{-1} } \geq 0 \ and \ \gamma {}_f = - \gamma {}_{f^{-1} } \in [ - \pi /2 , \pi /2] . 
\] 
Let $ \rho $ be a unitary representation of $ \Gamma $ 
and $d$ its degree. 
The {\em $L$-function} of $G$ associated with 
$ \rho $ and $ \alpha $ is defined by 
\[
{\bf Z}_H (G, \lambda , \rho , \alpha )= \det ( {\bf I}_{2m d } 
- \sum_{h \in \Gamma } {\rho} (h) \bigotimes {\bf U} {}_h )^{-1} . 
\]

If $\rho ={\bf 1}$ is the identity representation of $\Gamma $, then 
\[
{\bf Z}_H (G, \lambda , {\bf 1} , \alpha )= \det ( {\bf I}_{2m } - {\bf U} )^{-1} . 
\] 

A determinant expression for the $L$-function of $G$ associated with 
$ \rho $ and $ \alpha $ is given as follows. 
For $1 \leq i,j \leq n$, the {\em $(i,j)$-block} ${\bf F}_{i,j} $ of a 
$dn \times dn$ matrix ${\bf F}$ is the submatrix of ${\bf K}$ 
consisting of $d(i-1)+1, \ldots , di$ rows and 
$d(j-1)+1, \ldots , dj$ columns.

\begin{theorem}
Let $G$ be a connected graph with $n$ vertices and $m$ unoriented edges, 
$ \Gamma $ be a finite group and 
$ \alpha : D(G) \longrightarrow \Gamma $ be an ordinary voltage assignment. 
If $ \rho $ is a unitary representation of $ \Gamma $ and $d$ is 
the degree of $ \rho $, then the reciprocal of the $L$-function of $G$ 
associated with $ \rho $ and $ \alpha $ is
\[
{\bf Z}_H (G, \lambda , \rho , \alpha )^{-1} 
=\frac{ 2^{md} (-1)^{nd} }{\prod_{u \in V(G)} ( H_{uu} - \lambda -i \Gamma {}_u )^d } 
\det ( \lambda {\bf I}_{np} - \sum_{g \in \Gamma } \rho (g) \bigotimes {\bf H}_g - 
{\bf I}_d \bigotimes {\rm diag} ( {\bf H} )) . 
\] 
\end{theorem}

{\bf Proof}.  The argument is an analogue of Watanabe and Fukumizu's method \cite{WF}. 

Let $V(G)= \{ v_1 , \ldots , v_n \} $ and $D(G)= \{ e_1, \ldots , e_m , e_{m+1},$ 
$\ldots , e_{2m} \} $ such that 
$ e_{m+i} = e^{-1}_i (1 \leq i \leq m)$. 
Note that the $(e,f)$-block 
$(\sum_{g \in \Gamma } {\bf U}_g \bigotimes \rho (g) )_{ef} $ 
of $ \sum_{g \in \Gamma } {\bf U}_g \bigotimes \rho (g)$ is given by 
\[ 
(\sum_{g \in \Gamma } {\bf U}_g \bigotimes \rho (g) )_{ef} 
=\left\{
\begin{array}{ll}
\rho ( \alpha (e)) \sigma {}^{(t(e))}_{ef} & \mbox{if $t(e)=o(f)$, } \\
{\bf 0}_d  & \mbox{otherwise.}
\end{array}
\right.
\] 

For $g \in \Gamma $, two $2m \times 2m$ matrices 
${\bf B}_g =( {\bf B}^{(g)}_{ef} )_{e,f \in D(G)} $ and 
${\bf J}_g =( {\bf J}^{(g)}_{ef} )_{e,f \in D(G)} $ 
are defined as follows: 
\[
{\bf B}^{(g)}_{ef} =\left\{
\begin{array}{ll}
x_{o(f)} w(e) w(f) & \mbox{if $t(e)=o(f)$ and $ \alpha (e)=g$, } \\
0 & \mbox{otherwise, }
\end{array}
\right.
\  
{\bf J}^{(g)}_{ef} =\left\{
\begin{array}{ll}
1 & \mbox{if $f= e^{-1} $ and $\alpha (e)=g$, } \\
0 & \mbox{otherwise.}
\end{array}
\right.
\]
Then we have   
\[
{\bf U}_g =i {\bf J}_g -{\bf B}_g \  for \  g \in \Gamma .  
\]

Let ${\bf K} =( {\bf K}_{ij} )$ ${}_{1 \leq i \leq 2m; 
1 \leq j \leq n} $ be the $2md \times nd$ matrix defined 
as follows: 
\[
{\bf K}_{ij} :=\left\{
\begin{array}{ll}
x_{v_j} w(e_i) {\bf I}_d & \mbox{if $o( e_i )= v_j $, } \\
{\bf 0}_d & \mbox{otherwise. } 
\end{array}
\right.
\]
Furthermore, we define two $2md \times nd$ matrices 
${\bf L} =( {\bf L}_{ij} )_{1 \leq i \leq 2m; 
1 \leq j \leq n} $ and 
${\bf M} =( {\bf M}_{ij} )_{1 \leq i \leq 2m; 
1 \leq j \leq n} $ as follows: 
\[
{\bf L}_{ij} :=\left\{
\begin{array}{ll}
w( e_j ) \rho ( \alpha (e_i )) & \mbox{if $t( e_i )= v_j $, } \\
{\bf 0}_d & \mbox{otherwise, } 
\end{array}
\right.
\ 
{\bf M}_{ij} :=\left\{
\begin{array}{ll}
w(e_i ) {\bf I}_d & \mbox{if $o( e_i )= v_j $, } \\
{\bf 0}_d & \mbox{otherwise. } 
\end{array}
\right.
\]
Then we have 
\[
{\bf K} = {\bf M} ( {\bf X} \bigotimes {\bf I}_d )= {\bf M} {\bf X}_d , 
\]
where 
\[
{\bf X}_d = {\bf X} \bigotimes {\bf I}_d . 
\]
Furthermore, we have   
\begin{equation}
{\bf L} \ {}^t {\bf K} = \sum_{h \in \Gamma } {\bf B} {}_h \bigotimes \rho (h)= {\bf B}_{\rho }     
\end{equation}
and 
\begin{equation}
{}^t {\bf M} {\bf L} = \sum_{g \in \Gamma } {\bf H}_g \bigotimes \rho (g) ,  
\end{equation} 
where 
\[
{\bf B}_{\rho } = \sum_{g \in \Gamma } {\bf B}_g \bigotimes \rho (g) . 
\]
Thus, 
\[
\begin{array}{rcl}
\  &   & \det ( {\bf I}_{2md} - u \sum_{g \in \Gamma } \rho (g) \bigotimes {\bf U}_g )
= \det ( {\bf I}_{2md} - u \sum_{g \in \Gamma } {\bf U}_g \bigotimes \rho (g)) \\ 
\  &   &                \\ 
\  & = & \det ( {\bf I}_{2md} - u \sum_{g \in \Gamma } (i {\bf J}_g -{\bf B}_g ) \bigotimes \rho (g)) \\ 
\  &   &                \\ 
\  & = & \det ( {\bf I}_{2md} -i u \sum_{g \in \Gamma } {\bf J}_g \bigotimes \rho (g) 
+u \sum_{g \in \Gamma } {\bf B}_g \bigotimes \rho (g)) . 
\end{array}
\]

Now, let 
\[
{\bf J}_{\rho } = \sum_{g \in \Gamma } {\bf J}_g \bigotimes \rho (g) . 
\]
Note that 
\[
{\bf J}^2_{\rho } = {\bf I}_{2md} . 
\]
Then we have 
\[
\begin{array}{rcl}
\  &   & \det ( {\bf I}_{2md} - u \sum_{g \in \Gamma } \rho (g) \bigotimes {\bf U}_g ) \\ 
\  &   &                \\ 
\  & = & \det ( {\bf I}_{2md} -i u {\bf J}_{ \rho } +u {\bf B}_{\rho } ) \\ 
\  &   &                \\ 
\  & = & \det ( {\bf I}_{2md} +u {\bf B}_{\rho } ( {\bf I}_{2md} -i u {\bf J}_{ \rho })^{-1} ) 
\det ( {\bf I}_{2md} -i u {\bf J}_{ \rho } ) \\ 
\  &   &                \\ 
\  & = & \det ( {\bf I}_{2md} +u {\bf L} \ {}^t {\bf K} ( {\bf I}_{2md} -i u {\bf J}_{ \rho })^{-1} ) 
\det ( {\bf I}_{2md} -i u {\bf J}_{ \rho } ) \\ 
\  &   &                \\ 
\  & = & \det ( {\bf I}_{nd} +u \ {}^t {\bf K} ( {\bf I}_{2md} -i u {\bf J}_{ \rho })^{-1} {\bf L} ) 
\det ( {\bf I}_{2md} -i u {\bf J}_{ \rho } ) . 
\end{array}
\]

But, we have 
\[
\det ( {\bf I}_{2md} -i u {\bf J}_{ \rho } )= \det ( 
\left[
\begin{array}{cccc}
{\bf I}_d & -iu \rho ( \alpha ( e_1 )) &  & {\bf 0} \\ 
-iu \rho ( \alpha ( e^{-1}_1 )) & {\bf I}_{d} & & \\ 
  &  & \ddots & \\
{\bf 0} &   &   &  
\end{array}
\right]
)
=(1+u^2 )^{md} . 
\]
Furthermore, we have 
\[
\begin{array}{rcl}
\  &   & ( {\bf I}_{2md} -i u {\bf J}_{ \rho })^{-1} \\ 
\  &   &                \\ 
\  & = & 
\left[
\begin{array}{cccc}
{\bf I}_d & -iu \rho ( \alpha ( e_1 )) &  & {\bf 0} \\ 
-iu \rho ( \alpha ( e^{-1}_1 )) & {\bf I}_{d} & & \\ 
  &  & \ddots & \\
{\bf 0} &   &   &  
\end{array}
\right]^{-1} \\
\  &   &                \\ 
\  & = & \frac{1}{1+ u^2 } 
\left[
\begin{array}{cccc}
{\bf I}_d & iu \rho ( \alpha ( e_1 )) &  & {\bf 0} \\ 
iu \rho ( \alpha ( e^{-1}_1 )) & {\bf I}_{d} & & \\ 
  &  & \ddots & \\
{\bf 0} &   &   &  
\end{array}
\right]
\\ 
\  &   &                \\ 
\  & = & \frac{1}{1+ u^2 } ( {\bf I}_{2md} +i u {\bf J}_{ \rho } ) . 
\end{array}
\]
Thus, we have 
\[
\begin{array}{rcl}
\  &   & \det ( {\bf I}_{2md} - u \sum_{g \in \Gamma } \rho (g) \bigotimes {\bf U}_g ) \\ 
\  &   &                \\ 
\  & = & (1+ u^2 )^{md} \det ( {\bf I}_{nd} +u/(1+ u^2 ) \ {}^t {\bf K} 
( {\bf I}_{2md} +i u {\bf J}_{ \rho }) {\bf L} ) \\ 
\  &   &                \\ 
\  & = & (1+ u^2 )^{md-nd} \det ((1+ u^2 ) {\bf I}_{nd} +u \ {}^t {\bf K} {\bf L} 
+i u^2 \ {}^t {\bf K} {\bf J}_{ \rho } {\bf L} ) . 
\end{array}
\]

Now, we have 
\[
{}^t {\bf K} {\bf L} = {\bf X}_d \ {}^t {\bf M} {\bf L} 
={\bf X}_d \sum_{g \in \Gamma } {\bf H}_g \bigotimes \rho (g) . 
\]
Furthermore, 
\[
{}^t {\bf K} {\bf J}_{ \rho } {\bf L} = {\bf X}_d \ {}^t {\bf M} {\bf J}_{ \rho } {\bf L} . 
\]
Then we have 
\[
\begin{array}{rcl}
\  &   & ( {}^t {\bf M} {\bf J}_{ \rho } {\bf L})_{uv} \\ 
\  &   &                \\ 
\  & = & \delta {}_{uv} \sum_{o(e)=u} ( {}^t {\bf M} )_{ue} ({\bf J}_{ \rho } )_{e e^{-1} } ( {\bf L} )_{e^{-1} v} \\ 
\  &   &                \\ 
\  & = & \delta {}_{uv} \sum_{o(e)=u} w(e) {\bf I}_d \rho ( \alpha (e)) w( e^{-1} ) \rho ( \alpha ( e^{-1} )) \\
\  &   &                \\ 
\  & = & \delta {}_{uv} \sum_{o(e)=u} \sqrt{h_e } e^{i \gamma {}_e } \sqrt{h_e } e^{-i \gamma {}_e } {\bf I}_d \\ 
\  &   &                \\ 
\  & = &  \delta {}_{uv} \sum_{o(e)=u} h_e {\bf I}_d = \delta {}_{uv} \Gamma {}_u {\bf I}_d . 
\end{array}
\]
Thus, 
\[
{}^t {\bf K} {\bf J}_{ \rho } {\bf L} = {\bf X} ( {\bf D}_{\Gamma } \bigotimes {\bf I}_d ) , 
\]
where 
\[
{\bf D}_{\Gamma } = 
\left[
\begin{array}{ccc}
\Gamma {}_{v_1 } &  & 0 \\ 
  & \ddots & \\
0 &   & \Gamma {}_{v_n }   
\end{array}
\right]
. 
\]
Therefore, it follows that 
\[
\begin{array}{rcl}
\  &   & \det ( {\bf I}_{2md} - u \sum_{g \in \Gamma } \rho (g) \bigotimes {\bf U}_g ) \\ 
\  &   &                \\ 
\  & = & (1+ u^2 )^{(m-n)d} \det ((1+ u^2 ) {\bf I}_{nd} +u {\bf X}_d \sum_{g \in \Gamma } {\bf H}_g \bigotimes \rho (g) 
+i u^2 {\bf X}_d ( {\bf D}_{\Gamma } \bigotimes {\bf I}_d )) . 
\end{array}
\]

Substituting $u=1$, we obtain 
\[
\begin{array}{rcl}
\  &   & \det ( {\bf I}_{2md} - \sum_{g \in \Gamma } \rho (g) \bigotimes {\bf U}_g ) \\ 
\  &   &                \\ 
\  & = & 2^{(m-n)d} \det (2 {\bf I}_{nd} + {\bf X}_d \sum_{g \in \Gamma } {\bf H}_g \bigotimes \rho (g) 
+i {\bf X}_d ( {\bf D}_{\Gamma } \bigotimes {\bf I}_d )) \\ 
\  &   &                \\ 
\  & = & 2^{(m-n)d} \det ({\bf X}_d ) \det (2 {\bf X}^{-1}_d + \sum_{g \in \Gamma } {\bf H}_g \bigotimes \rho (g) 
+i {\bf D}_{\Gamma } \bigotimes {\bf I}_d ) .  
\end{array}
\]
Then we have 
\[
\det ({\bf X}_d )= \det ({\bf X} \bigotimes {\bf I}_d ) =(\det ({\bf X} ))^d 
= \frac{2^{nd}}{ \prod_{u \in V(G)} ( H_{uu} - \lambda -i \Gamma {}_u )^d } . 
\]
Furthermore, since 
\[
{\bf X}^{-1}_d = {\bf X}^{-1} \bigotimes {\bf I}_d , 
\]
we have 
\[
\begin{array}{rcl}
(2 {\bf X}^{-1}_d +i {\bf D}_{\Gamma } \bigotimes {\bf I}_d )_{uu} & = & 
(2 \frac{H_{uu} - \lambda -i \Gamma {}_u }{2} +i \Gamma {}_u ) \bigotimes {\bf I}_d \\
\  &   &                \\ 
\  & = & (H_{uu} - \lambda ) \bigotimes {\bf I}_d . 
\end{array}
\]
That is, 
\[
2 {\bf X}^{-1}_d +i {\bf D}_{\Gamma } \bigotimes {\bf I}_d 
=- \lambda {\bf I}_{nd} + {\rm diag} ( {\bf H} ) \bigotimes {\bf I}_d . 
\]
Therefore, it follows that 
\[ 
\begin{array}{rcl}
\  &   & \det ( {\bf I}_{2md} - \sum_{g \in \Gamma } \rho (g) \bigotimes {\bf U}_g ) \\ 
\  &   &                \\ 
\  & = & \frac{2^{md}}{ \prod_{u \in V(G)} ( H_{uu} - \lambda -i \Gamma {}_u )^d }  
\det (- \lambda {\bf I}_{nd} + \sum_{g \in \Gamma } {\bf H}_g \bigotimes \rho (g) + {\rm diag} ( {\bf H} ) \bigotimes {\bf I}_d ) \\ 
\  &   &                \\ 
\  & = & \frac{(-1)^{nd} 2^{md}}{ \prod_{u \in V(G)} ( H_{uu} - \lambda -i \Gamma {}_u )^d }  
\det ( \lambda {\bf I}_{nd} - \sum_{g \in \Gamma } \rho (g) \bigotimes {\bf H}_g 
- {\bf I}_d \bigotimes {\rm diag} ( {\bf H} ) ) . 
\end{array}
\]
$\Box$

By Theorems 6 and 7 the following result holds.

\newtheorem{corollary}{Corollary}
\begin{corollary}
Let $G$ be a connected graph with $m$ edges, $ \Gamma $ be a finite group and 
$ \alpha : D(G) \longrightarrow \Gamma $ be an ordinary voltage assignment. 
Then  
\[
\det ( {\bf I}_{2mp} - \tilde{{\bf U}} ( \lambda ))= 
\prod_{ \rho } {\bf Z}_H (G,\lambda , \rho , \alpha )^{- \deg \rho } , 
\]
where $ \rho $ runs over all inequivalent irreducible representations 
of $ \Gamma $ and $p= \mid \Gamma \mid $. 
\end{corollary}

\section{Example}

We give an example.
Let $G=K_3$ be the complete graph with three vertices $1,2,3$ and 
six arcs $e_1 ,e_2 , e_3 , e^{-1}_1 ,e^{-1}_2 , e^{-1}_3 $, 
where $e_1 =( v_1 ,v_2 ), e_2 = ( v_2 ,v_3 ), e_3=( v_3 ,v_1 )$.  
Furthermore, let 
\[
{\bf H} =
\left[ 
\begin{array}{ccc}
a & b e^{2i \alpha } & b e^{2i \alpha } \\
b e^{-2i \alpha } & a & b e^{2i \alpha } \\ 
b e^{-2i \alpha } & b e^{-2i \alpha } & a 
\end{array} 
\right] 
,  
\]
where $a>0$, $b>0$ and $ \alpha \in [- \frac{\pi}{2} , \frac{\pi }{2} )$. 
Then we have 
\[
x_{1} = x_{2} =x_{3} = \frac{2}{a- \lambda -2ib} . 
\]
Set $x= \frac{2}{a- \lambda -2ib} $. 
Considering ${\bf U} ( \lambda ) $ under the order 
$e_1 ,e_2 , e_3 , e^{-1}_1 ,e^{-1}_2 , e^{-1}_3 $, 
we have 
\[
{\bf U} ( \lambda ) =
\left[ 
\begin{array}{cccccc}
-xb e^{2i \alpha } & -xb e^{2i \alpha } & -xb e^{2i \alpha } & i-xb & -xb & -xb \\
-xb e^{2i \alpha } & -xb e^{2i \alpha } & -xb e^{2i \alpha } & -xb & i-xb & -xb \\
-xb e^{2i \alpha } & -xb e^{2i \alpha } & -xb e^{2i \alpha } & -xb & -xb & i-xb \\
i-xb & -xb & -xb & -xb e^{-2i \alpha } & -xb e^{-2i \alpha } & -xb e^{-2i \alpha } \\
-xb & i-xb & -xb & -xb e^{-2i \alpha } & -xb e^{-2i \alpha } & -xb e^{-2i \alpha } \\
-xb & -xb & i-xb & -xb e^{-2i \alpha } & -xb e^{-2i \alpha } & -xb e^{-2i \alpha }  
\end{array} 
\right] 
. 
\]
By Theorem 4, we have 
\[
\begin{array}{rcl}
\det ({\bf I}_6 - {\bf U} ( \lambda ) ) & = & \frac{2^3 (-1)^3}{(a- \lambda -2ib )^3} 
\det ( \lambda {\bf I}_3 - {\bf H} ) \\
\  &   &                \\ 
\  & = & \frac{-8}{(a- \lambda -2ib )^3} 
\left[ 
\begin{array}{ccc}
\lambda -a & -b e^{2i \alpha } & -b e^{2i \alpha } \\
-b e^{-2i \alpha } & \lambda -a & -b e^{2i \alpha } \\ 
-b e^{-2i \alpha } & -b e^{-2i \alpha } & \lambda -a 
\end{array} 
\right] 
\\ 
\  &   &                \\ 
\  & = & \frac{-8}{(a- \lambda -2ib )^3} 
\{ (\lambda -a)^3 -3b^2 ( \lambda -a)-b^3 ( e^{2i \alpha } +e^{-2i \alpha } ) \} \\
\  &   &                \\ 
\  & = & \frac{-8}{(a- \lambda -2ib )^3} 
\{ (\lambda -a)^3 -3b^2 ( \lambda -a)-2b^3 \cos 2 \alpha \} . 
\end{array}
\]

Next. let $\Gamma = \mathbb{Z}_3 = \{ 1, \tau , { \tau }^2 \} ( { \tau }^3 =1)$ be 
the cyclic group of order 3, and let 
$ \alpha : D(K_3 ) \longrightarrow \mathbb{Z}_3$ be the ordinary voltage 
assignment such that $ \alpha ( e_1 )= \tau $, $ \alpha ( e^{-1}_1 )= \tau {}^2 $ 
and $ \alpha ( e_2 )= \alpha ( e^{-1}_2 )= \alpha ( e_3 )= \alpha (e^{-1}_3 )=1$. 
Then the $\mathbb{Z}_3 $-coverng $K^{\alpha }_3 $ of $K_3 $ is the cycle graph of length 9. 

The characters of $\mathbb{Z}_3$ are given as follows: 
$\chi {}_i ( \tau {}^j )=( \xi {}^i )^j $, $ 0 \leq i, j \leq 2$, 
where $ \xi = \frac{-1+ \sqrt{-3} }{2} $. 
Then  we have 
\[
{\bf H}_1 =
\left[ 
\begin{array}{ccc}
0 & 0 & b e^{2i \alpha } \\
 & 0 & b e^{2i \alpha } \\ 
b e^{-2i \alpha } & b e^{-2i \alpha } & 0 
\end{array} 
\right] 
, 
{\bf H}_{\tau } =
\left[ 
\begin{array}{ccc}
0 & b e^{2i \alpha } & 0 \\
0 & 0 & 0 \\ 
0 & 0 & 0 
\end{array} 
\right] 
, 
{\bf H}_{\tau {}^2 } =
\left[ 
\begin{array}{ccc}
0 & 0 & 0 \\
b e^{-2i \alpha } & 0 & 0 \\ 
0 & 0 & 0 
\end{array} 
\right] 
.
\]

Now, by Theorem 7, 
\[
\begin{array}{rcl}
\  &  & {\zeta}_H (K_3, \lambda , \chi {}_1 , \alpha )^{-1} =  
\frac{2^3 (-1)^3}{(a- \lambda -2ib )^3} 
\det ( \lambda {\bf I}_3 - \sum^2_{j=0} \chi {}_1 ( \tau {}^j ) {\bf H}_{ \tau {}^j } 
- {\rm diag} ( {\bf H} )) \\
\  &   &                \\ 
\  & = & \frac{-8}{(a- \lambda -2ib )^3} 
\left[ 
\begin{array}{ccc}
\lambda -a & -b \xi e^{2i \alpha } & -b e^{2i \alpha } \\
-b \xi {}^2 e^{-2i \alpha } & \lambda -a & -b e^{2i \alpha } \\ 
-b e^{-2i \alpha } & -b e^{-2i \alpha } & \lambda -a 
\end{array} 
\right] \\
\  &   &                \\ 
\  & = & \frac{-8}{(a- \lambda -2ib )^3} 
\{ (\lambda -a)^3 -3b^2 ( \lambda -a)-b^3 ( \xi e^{2i \alpha } + \xi {}^2 e^{-2i \alpha })  \} \\
\  &   &                \\ 
\  & = & \frac{-8}{(a- \lambda -2ib )^3} 
\{ (\lambda -a)^3 -3b^2 ( \lambda -a)-2b^3 \cos 2 ( \alpha +\pi /3) \} . 
\end{array}
\] 
Similarly, we have 
\[
\begin{array}{rcl}
\  &  & {\zeta}_H (K_3, \lambda , \chi {}_2 , \alpha )^{-1} =  
\frac{2^3 (-1)^3}{(a- \lambda -2ib )^3} 
\det ( \lambda {\bf I}_3 - \sum^2_{j=0} \chi {}_2 ( \tau {}^j ) {\bf H}_{ \tau {}^j } 
- {\rm diag} ( {\bf H} )) \\
\  &   &                \\ 
\  & = & \frac{-8}{(a- \lambda -2ib )^3} 
\left[ 
\begin{array}{ccc}
\lambda -a & -b \xi {}^2 e^{2i \alpha } & -b e^{2i \alpha } \\
-b \xi e^{-2i \alpha } & \lambda -a & -b e^{2i \alpha } \\ 
-b e^{-2i \alpha } & -b e^{-2i \alpha } & \lambda -a 
\end{array} 
\right] \\
\  &   &                \\ 
\  & = & \frac{-8}{(a- \lambda -2ib )^3} 
\{ (\lambda -a)^3 -3b^2 ( \lambda -a)-b^3 ( \xi {}^2 e^{2i \alpha } + \xi e^{-2i \alpha } ) \} \\
\  &   &                \\ 
\  & = & \frac{-8}{(a- \lambda -2ib )^3} 
\{ (\lambda -a)^3 -3b^2 ( \lambda -a)-2b^3 \cos 2 ( \alpha +2 \pi /3) \} . 
\end{array}
\]

By Corollary 1, it follows that 
\[
\begin{array}{rcl}
\  &  & \det ({\bf I}_{18} - \tilde{{\bf U}} ( \lambda ) )= 
\det ({\bf I}_6 - {\bf U} ( \lambda ) ) \zeta {}_H (K_3,\lambda , \chi {}_1, \alpha )^{-1}  
{\zeta}_S (K_H, \lambda , \chi {}_2 , \alpha )^{-1} \\ 
\  &   &                \\ 
\  & = & \frac{-512}{(a- \lambda -2ib )^9} \{ (\lambda -a)^3 -3b^2 ( \lambda -a)-2b^3 \cos 2 \alpha \} \\
\  &   &                \\ 
\  & \times & \{ (\lambda -a)^3 -3b^2 ( \lambda -a)-2b^3 \cos 2 ( \alpha +\pi /3) \} 
\{ (\lambda -a)^3 -3b^2 ( \lambda -a)-2b^3 \cos 2 ( \alpha +2 \pi /3) \} . 
\end{array}
\]

\end{document}